\documentclass[final]{siamltex}
\usepackage{amsfonts,graphicx}
\usepackage[leqno]{amsmath}

\def\gg{G^{(\gamma)}}

\def\ja1x{P_{n-1}^{(\alpha,\beta)}(x)}

\def\CC{\mathcal{C}}
\def\jacw8{W_{\alpha,\beta}}

\def\ie{{\it i.e.}}
\def\etal{{\it et al.}}
\def\f12{\frac{1}{2}}

\newcommand{\boxit}[1]{\vbox{\hrule\hbox{\strut \vrule \vbox{ \kern6pt #1}
\vrule}\hrule}}
\newtheorem{rmk}{Remark}

\begin{document}
\title{Gegenbauer tau Methods with and without Spurious
Eigenvalues}
\author{Marios Charalambides%
\thanks{Department of Business Administration, Frederick University Cyprus, 
7 Yianni Frederickou Street, Pallouriotissa, PO Box 24729, 1303
Nicosia,
Cyprus (\texttt{bus.chm@fit.ac.cy}).} 
\and Fabian Waleffe%
\thanks{Department of Mathematics, University of Wisconsin, Madison, WI 53706, USA (\texttt{waleffe@math.wisc.edu}). This work was supported in part by
NSF grant DMS-0204636.}}
\date{2007/09/19 MC}
\maketitle

\begin{abstract}
It is proven that a class of Gegenbauer tau approximations to a 4th
order differential eigenvalue problem of hydrodynamic type provide
real, negative, and distinct eigenvalues, as is the case for the exact solutions. 
This class of Gegenbauer tau methods includes Chebyshev and Legendre
Galerkin and `inviscid' Galerkin  but does not include Chebyshev and
Legendre tau. Rigorous and numerical results show that the results
are sharp: positive or complex eigenvalues arise outside of this
class. The widely used modified tau approach is proved to be
equivalent to the Galerkin method.
\end{abstract}



\pagestyle{myheadings} \thispagestyle{plain} \markboth{MARIOS
CHARALAMBIDES AND FABIAN WALEFFE}{GEGENBAUER TAU METHODS WITH AND
WITHOUT SPURIOUS EIGENVALUES}

\section{Introduction}

\label{sec1}
The Chebyshev tau method used by Orszag  \cite{O71} to obtain exponentially accurate solutions of the Orr-Sommerfeld equation yields two 
eigenvalues with large positive real parts.
Such 
eigenvalues also occur for  Stokes modes in a channel given by
the fourth order differential equation
\begin{equation}
(D^2-\alpha^2)^2 \, u = \lambda (D^2-\alpha^2) u  \label{Stokes}
\end{equation}
with the boundary conditions $u=Du=0$ at $x=\pm 1$,
where $\lambda$ is the eigenvalue, $u=u(x)$ is the eigenfunction, $D=d/dx$ and $\alpha$ is a real wavenumber. 
The Stokes eigenvalues $\lambda$ are real and negative as can be checked by multiplying  (\ref{Stokes}) by $u^*$, the complex conjugate of $u(x)$, and integrating by parts twice using the no-slip boundary conditions. In fact the Stokes spectrum has been known analytically since Rayleigh \cite[\S 26.1]{DR81}. Yet, the Chebyshev tau method applied to (\ref{Stokes}) yields 2 eigenvalues with large positive real parts, for any order of approximation and for any numerical accuracy. Such eigenvalues  are obviously \emph{spurious} for (\ref{Stokes}). 
 Gottlieb and Orszag \cite[Chap.\ 13]{GO77} introduced the eigenvalue problem
\begin{equation}
\begin{aligned} D^4 u = \lambda \, D^2 u \quad &\mbox{in} \; -1 \le x \le 1, \\
u=Du=0 \quad &\mbox{at} \quad x=\pm 1, \end{aligned}
\label{testevp}
\end{equation}
as an even simpler 1D model of incompressible fluid flow. This
is the $\alpha \to 0$ limit of the eigenvalue problem (\ref{Stokes}) and of the Orr-Sommerfeld equation \cite{O71}.
For any fixed $\alpha$, problem (\ref{testevp}) is also the asymptotic equation for large $\lambda$
solutions of the Stokes and Orr-Sommerfeld equations.
The eigensolutions of (\ref{testevp}) are known analytically. They
consist of even modes 
$\displaystyle u(x) = 1-\cos(n \pi x)/\cos(n \pi)$ with
$\lambda=-n^2 \pi^2$
and odd modes 
$\displaystyle u(x) = x - {\sin(q_n x)}/{\sin(q_n)}$, with
$\lambda=-q_n^2$,
where $q_n = \tan q_n$, so that $n \pi < q_n < (2n+1) \pi/2$,
$\forall$ integer $n>0$. The key properties of these
solutions are that the eigenvalues are \emph{real, negative, and distinct,} and the even and odd mode eigenvalues interlace. These properties
also hold for Stokes modes, the solutions of (\ref{Stokes}), but not for Orr-Sommerfeld modes.

The Chebyshev tau method provides spectrally accurate approximations
to the lower magnitude eigenvalues but it also yields two large
positive eigenvalues for problems (\ref{Stokes}) and (\ref{testevp}) \cite[Table 13.1]{GO77}.
Those positive eigenvalues are clearly \textit{spurious} since it is known that  (\ref{Stokes}) and (\ref{testevp}) should only have negative eigenvalues.
The Chebyshev tau method yields 2 spurious eigenvalues for no-slip
(\emph{a.k.a.} `clamped') boundary conditions $u=Du=0 $ at $x=\pm
1$, but none for the free-slip boundary conditions $u=D^2u=0$ at
$x=\pm 1$. The latter problem reduces to the 2nd order problem $D^2 v =
\lambda v$ with $v(\pm 1)=0$ for which a class of Jacobi and
Gegenbauer tau methods has been proven to yield real, negative and
distinct eigenvalues \cite{CCW,CW}. For mixed boundary conditions,
\textit{e.g.} $u(\pm1)=Du(-1)=D^2u(1)=0$ there is one spurious
eigenvalue (this is a numerical observation).

\begin{figure}
\centering
\includegraphics[height=1.95truein]{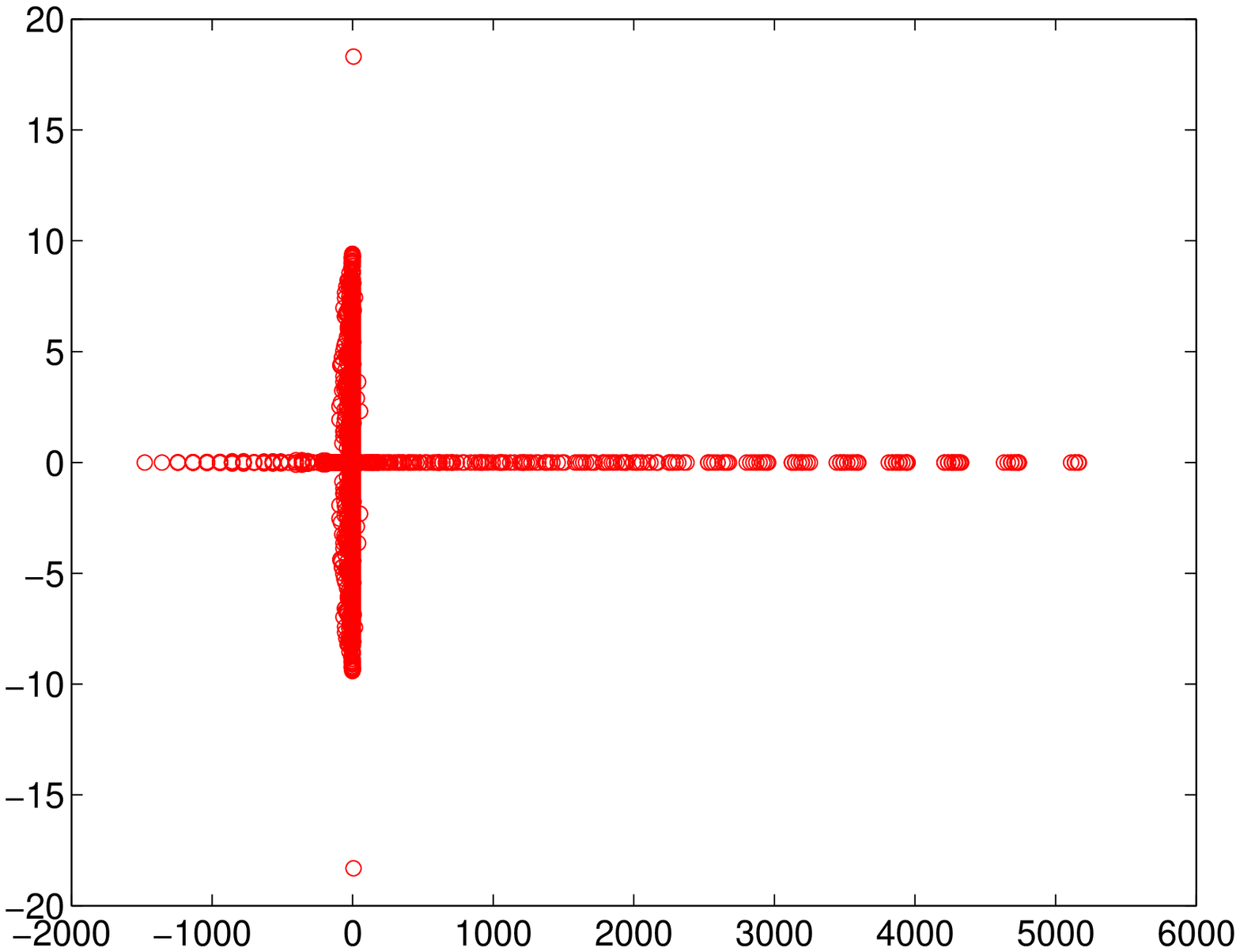}
\includegraphics[height=1.95truein]{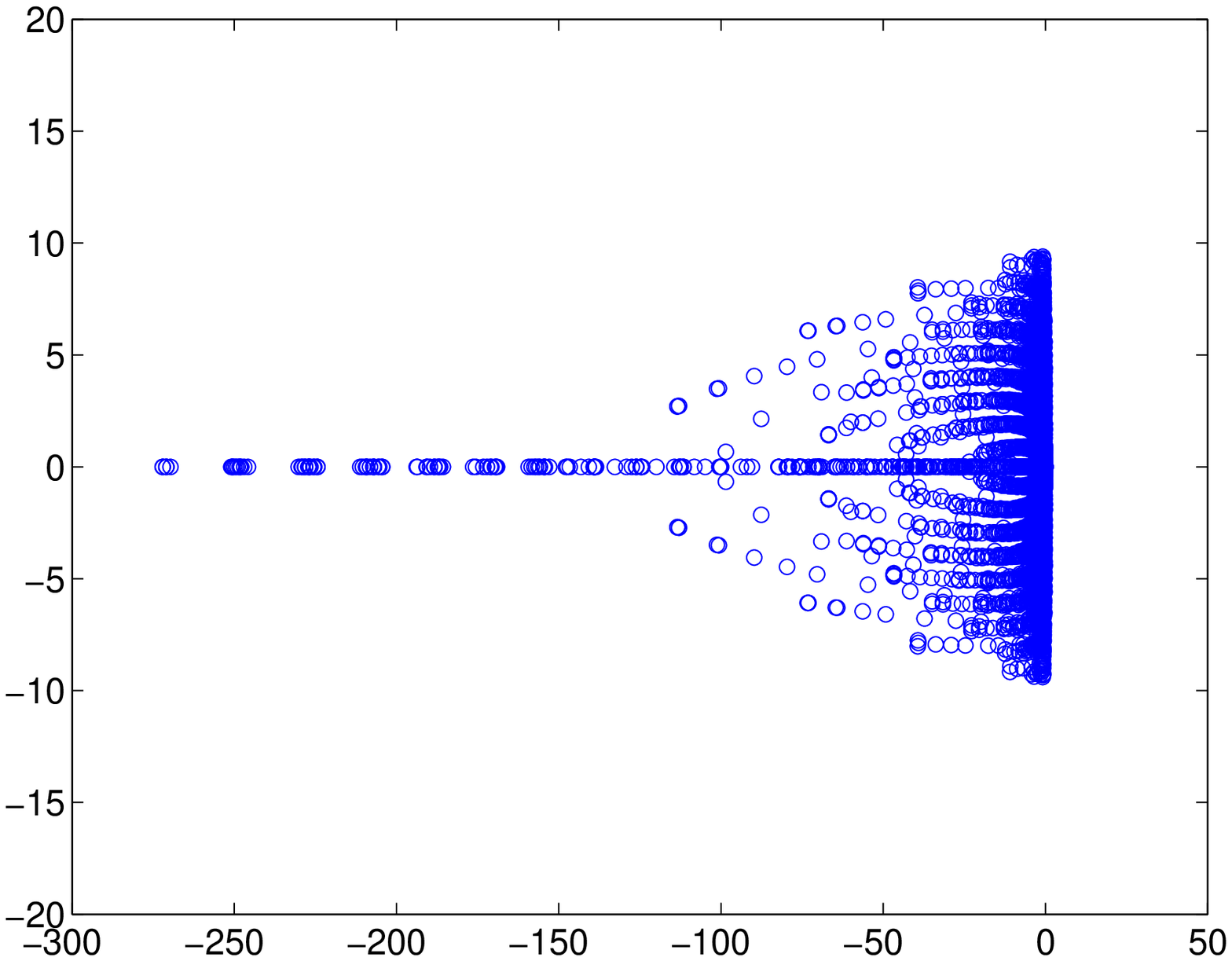}
\caption{Eigenvalues of a 3D steady state solution of the Navier-Stokes equations for plane Couette flow \cite{WGW} computed with Chebyshev tau (left) and Chebyshev Galerkin (right) for identical resolutions (8773 modes after symmetry reductions. The solutions themselves are indistinguishable). Note the difference in horizontal scales. Chebyshev tau produces 274 eigenvalues with positive real parts, 273 of which are spurious. Chebyshev Galerkin returns only one positive eigenvalue, the physical one, equal to 0.03681 at Reynolds number 1000 \cite[Fig.\ 4]{WGW}.}
\label{eigs-LBS}
\end{figure}

 For a 1D problem such as (\ref{Stokes}) or (\ref{testevp}), the spurious eigenvalues are easy to recognize and they appear
 as minor nuisances.  
 Boyd \cite[\S 7.6]{Boyd} even questions the value of distinguishing between `spurious' and numerically
inaccurate eigenvalues. 
However, in many applications, large \emph{negative} eigenvalues are inconsequential, while `spurious' \emph{positive} eigenvalues are very significant, and in higher dimensions, spurious
eigenvalues are not as easy to pick out and set aside. In a recent
application, 3D \emph{unstable} traveling wave solutions of the
Navier-Stokes equations were calculated with both free-slip and
no-slip boundary conditions, and anything in between,
  by Newton's method \cite{W01,W03}. In that application,
the Chebyshev tau method provides hundreds of spurious
unstable eigenvalues, depending on resolution and the exact type of
boundary conditions, not all of which have very large magnitudes (fig.\ \ref{eigs-LBS}, left). 
A simple change in the \emph{test} functions from Chebyshev
polynomials $T_n(x)$ to  $(1-x^2) T_n(x)$ or $(1-x^2)^2 T_n(x)$
eliminates all those spurious eigenvalues (fig.\ \ref{eigs-LBS}, right). This, and more, is
proven below for the test problem (\ref{testevp}) in the broader
context of Gegenbauer tau methods which include Chebyshev and
Legendre tau, and Chebyshev and Legendre Galerkin methods.

Another practical consequence of the spurious eigenvalues is that
the Chebyshev tau method is unconditionally unstable when applied to
the time-dependent version of (\ref{Stokes}) or (\ref{testevp}),
with $\partial/\partial t$ in place of $\lambda$. Such
time-dependent problems appear as building blocks in Navier-Stokes
simulations of channel-type flows. 
Gottlieb and Orszag \cite[p.\ 145]{GO77} proposed a \emph{modified tau} method for the
time-dependent problems and proved that the modified method was
stable for even solutions. 
The modified tau method (sect.\ \ref{numerical}) is a key idea behind several successful time integration schemes for the Navier-Stokes equations \cite[\S 7.3]{CHQZ}, \cite{KMM}.
The modified tau method amounts to using
2 more expansion polynomials for the 4th order differential operator
on the left hand sides of (\ref{testevp}) and (\ref{Stokes}) than
for the 2nd order operator on the right hand sides. 
 That modified tau method was adapted to eigenvalue problems by
Gardner, Trogdon and Douglass \cite{GTD89} and McFadden, Murray and
Boisvert \cite{MMB90}. McFadden \etal\  showed the equivalence
between the modified Chebyshev tau method and a Chebyshev Galerkin
method by direct calculation. Zebib \cite{Zebib87} had given numerical evidence that the Galerkin method removed spurious eigenvalues. The modified tau method idea was
adapted to the collocation formulation by Huang and Sloan \cite{HS94}. 

A heuristic `explanation' for spurious eigenvalues is that there is a `mismatch' between the number of boundary conditions applied to the 4th order operator on the left hand side of (\ref{testevp}) and those applied to the 2nd order operator on the right hand side. That interpretation fits with the modified tau method which uses two more polynomials for the 4th order operator than for the 2nd order operator. However it is incorrect since,  while the tau method for Chebyshev polynomials of the 1st kind $T_n(x)$ gives spurious eigenvalues, for instance, the tau method for Chebyshev polynomials of the 2nd kind $U_n(x)$ does not. 

All of these various methods are best seen in the context of the Gegenbauer class with residuals weighted by $W^{(\gamma)}(x) = (1-x^2)^{\gamma-1/2}$ (sect.\  \ref{residual}), where $\gamma=0$ corresponds to Chebyshev and $\gamma=1/2$ to Legendre polynomials.
Dawkins \etal \cite{DDD98} proved existence of spurious positive eigenvalues for (\ref{testevp}) when $\gamma < 1/2$.
The proof is straightforward. For (\ref{testevp}), the polynomial equation for $\mu=1/\lambda$ can be derived explicitly (sect.\ \ref{polys}). All coefficients of that polynomial are real and positive, except the constant term which is negative when $\gamma<1/2$. Hence there is one real positive $\mu$, and a `spurious' positive eigenvalue, when $\gamma<1/2$ (details are given in sect.\  \ref{zeroseven}). For $\gamma=1/2$,  the Legendre tau case, the constant term is zero, hence there is one $\mu=0$ eigenvalue, or a $\lambda=1/\mu=\infty$ eigenvalue.
Perturbation analysis shows that the $\lambda=\infty$ eigenvalues become very large positive eigenvalues for $\gamma<1/2$ and very large negative eigenvalues for $\gamma >1/2$.
We provide a quicker derivation of those results in section
\ref{Legendre}. Dawkins \etal's results do not prove that there are
no spurious eigenvalues for $\gamma>1/2$ since there could be
complex eigenvalues with positive real parts.
In section \ref{Zeros}, we prove that the Gegenbauer tau method applied to (\ref{testevp}) provides eigenvalues that are \emph{real, negative, and distinct} when $1/2 < \gamma \le 7/2$. This provides a complete characterization of the Gegenbauer tau spectrum for problem (\ref{testevp}).
Numerical calculations confirm that the range $1/2 < \gamma \le 7/2$ is sharp. Spurious positive
eigenvalues exist for $\gamma<1/2$ \cite{DDD98} and complex eigenvalues arise for
$\gamma > 7/2$, for sufficiently high polynomial order.
In section \ref{numerical}, we prove  that the modified tau method is mathematically equivalent to the Galerkin approach.

Obviously, $\gamma=1/2$ is a critical value for the weight function $(1-x^2)^{\gamma-1/2}$. The boundaries $x=\pm 1$ have infinite weight for $\gamma < 1/2$, and zero weight for $\gamma > 1/2$, but we do not know a valid heuristic explanation for `spurious' eigenvalues beyond that observation, if one exists. 
Section \ref{Legendre} provides further insights into the nature of the spurious eigenvalues and gives \emph{some} support for the view that `spurious' and numerically inaccurate eigenvalues are related. 
In figure \ref{eigs-LBS}, \emph{for fixed resolution}, the 273 spurious eigenvalues for $\gamma=0$ (Chebyshev tau), escape to $ + \infty$ as $\gamma \nearrow 1/2$ (Legendre tau). They come back from $-\infty$ as $\gamma$ increases beyond $1/2$. 
Thus, there is indeed a connection between large positive eigenvalues and large negative eigenvalues,
but whether we have `spurious' positive eigenvalues or inconsequential very negative eigenvalues is sharply controlled by $\gamma$, irrespective of the order of approximation $n$ (see also eqn.\ (\ref{spuriousleg})).


\section{Tau and Galerkin methods}

\label{residual}

\begin{definition}
A  Gegenbauer tau method approximates the solution $u(x)$ of a differential equation in  $-1 \le x \le 1$ by a polynomial of degree $n$, $u_n(x)$,  that satisfies the $m$ boundary conditions exactly. The remaining $n+1-m$ polynomial coefficients are determined by imposing that the residual be orthogonal to all polynomials of degree $n-m$ (or less) with respect to the Gegenbauer weight $W^{(\gamma)}(x) =(1-x^2)^{\gamma-1/2}$, with $ \gamma > -1/2$.
\end{definition}

For problem (\ref{testevp}), the residual 
\begin{equation}
R_{n-2}(x) \equiv \lambda D^2 u_n(x)-D^4u_n(x)
\label{Rn}
\end{equation}
is a polynomial of degree $n-2$ in $x$. The polynomial approximation $u_n(x)$ is determined from the 4 boundary conditions 
 $u_n(\pm 1) =Du_n(\pm 1)= 0$ and the requirement that $R_{n-2}(x)$ is orthogonal to \emph{all} polynomials $q_{n-4}(x)$ of
degree $n - 4$ or less with respect to the weight function $W^{(\gamma)}(x)=(1-x^2)^{\gamma-1/2} \ge
0$ in the interval $(-1, 1)$
\begin{equation}
\int_{-1}^1 R_{n-2}(x) \, q_{n-4}(x) \, W^{(\gamma)}(x) dx = 0, \quad \forall \, q_{n-4}(x). 
 \label{jactau}
\end{equation}
This provides $n-3$ equations  which together with the 4 boundary conditions yield
$n+1$ equations for the $n+1$ undetermined coefficients in the polynomial approximation
$u_n(x)$. For the Gegenbauer weight function $W^{(\gamma)}(x) =
(1-x^2)^{\gamma - 1/2}$, the residual can be written explicitly as 
\begin{equation}  \label{GTR}
R_{n-2}(x)=\tau_0 \lambda \, G_{n-2}^{(\gamma)}(x) +\tau_1 \lambda
\, G_{n-3}^{(\gamma)}(x)
\end{equation}
for some $x$-independent coefficients $\tau_0$ and $\tau_1$, where
$G_n^{(\gamma)}(x) $ is the Gegenbauer polynomial of degree $n$.
This follows from orthogonality of the Gegenbauer polynomials in
$-1 < x < 1$ with respect to the weight
$(1-x^2)^{\gamma - 1/2}$ which implies orthogonality of the
Gegenbauer polynomial of degree $k$ to \textit{any} polynomial of
degree $k-1$ or less with respect to that weight function.

Gegenbauer (\textit{a.k.a.} ultraspherical)  polynomials are a special subclass of the Jacobi polynomials \cite{AS}. The latter are the most general class of polynomial solutions of a
Sturm-Liouville eigenproblem that is singular at $\pm 1$, as
required for faster than algebraic convergence \cite{CHQZ}. Gegenbauer polynomials are the most general class of polynomials with the odd-even symmetry $\gg_{n}(x)= (-1)^n \gg_n(-x)$. This is a one parameter family of polynomials, with the parameter $\gamma > -1/2$.
Chebyshev polynomials correspond to $\gamma=0$ and Legendre polynomials to $\gamma=1/2$.
We use a (slightly) non-standard normalization of Gegenbauer
polynomials since the standard normalization \cite{AS} is singular in the Chebyshev case.
Some key properties of Gegenbauer polynomials used in this paper are given in Appendix \ref{GegAppend}. Note that if $\lambda = 0$ then, from (\ref{Rn}), the residual
$R_{n-2}(x)$ must be a polynomial of degree $n-4$ implying that
$\tau_0=\tau_1=0$ in (\ref{GTR}) and $D^4 u_n(x)=0$ for all $x$ in
$(-1,1)$.
The boundary conditions $u_n(\pm 1) =Du_n(\pm 1)= 0$ then imply
that $u_n(x)=0$ for all $x$ in $[-1,1]$, the trivial solution. Hence we can assume that $\lambda \ne 0$ in the Gegenbauer tau method applied to (\ref{testevp}).

Following common usage \cite{GO77,CHQZ,Zebib87,MMB90}, we have 
\begin{definition}
A Gegenbauer Galerkin method approximates the solution $u(x)$ of a differential equation in  $-1 \le x \le 1$ by a polynomial of degree $n$, $u_n(x)$,  that satisfies the $m$ boundary conditions exactly. The remaining $n+1-m$ polynomial coefficients are determined by imposing that the residual be orthogonal to all polynomials of degree $n$ (or less) that satisfy the homogeneous boundary conditions, 
with respect to the Gegenbauer weight $W^{(\gamma)}(x) =(1-x^2)^{\gamma-1/2}$, with $ \gamma > -1/2$.
\end{definition}

Strictly speaking, this a \emph{Petrov-Galerkin} method since the test functions are not identical to the trial functions because of the Gegenbauer weight $(1-x^2)^{\gamma-1/2}$ \cite{CHQZ}.

For problem (\ref{testevp}), $u_n(x)$ is determined from the
boundary conditions $u_n(\pm 1) =Du_n(\pm 1)= 0$ and orthogonality, with respect to weight $W^{(\gamma)}(x) = (1-x^2)^{\gamma-1/2}$,  of the residual (\ref{Rn})
to all polynomials of degree $n$ that vanish together with their
derivative at $x=\pm 1$.
Such polynomials can be written as $(1-x^2)^2 q_{n-4}(x) $ where
$q_{n-4}(x)$ is an arbitrary polynomial of degree $n-4$, and the weighted residual equations read
\begin{equation}
\int_{-1}^1 R_{n-2}(x) \, (1-x^2)^2 q_{n-4}(x) \, W^{(\gamma)}(x) dx = 0,  \quad \forall \, q_{n-4}(x).
\label{GG}
\end{equation}
The Gegenbauer Galerkin method is therefore equivalent to the tau
method for the weight $W^{(\gamma+2)}(x) = (1-x^2)^2 W^{(\gamma)}(x) $ and its residual has the explicit form
\begin{equation}  \label{GGR}
R_{n-2}(x)=\tau_0 \lambda \, G_{n-2}^{(\gamma+2)}(x) +\tau_1 \lambda \, G_{n-3}^{(\gamma+2)}(x).
\end{equation}
So a \emph{Chebyshev (or Legendre) Galerkin} method for clamped boundary conditions, $u_n(\pm 1)=Du_n(\pm 1)=0$,  is in fact a tau method for
Chebyshev (or Legendre) polynomials of the \emph{3rd kind} (proportional to the 2nd derivative of Chebyshev (or Legendre) polynomials (\ref{gege2})).
Since we consider a range of the Gegenbauer parameter $\gamma$,
the Gegenbauer tau method also includes some Gegenbauer Galerkin methods.

This suggests an intermediate method where the test functions are polynomials that vanish at $x=\pm 1$ (inviscid boundary conditions only).
\begin{definition}
The Gegenbauer `inviscid Galerkin' method determines $u_n(x)$ from the $4$ boundary conditions $u_n(\pm 1)=Du_n(\pm 1)=0$ and orthogonality of the residual  to all polynomials of degree $n-2$ that vanish at $x=\pm 1$, with respect to the weight function $W^{(\gamma)}(x)=(1-x^2)^{\gamma-1/2}$.
\end{definition}

Such test polynomials can be written in the form $(1-x^2) q_{n-4}(x)$ where
$q_{n-4}(x)$ is an arbitrary polynomial of degree $n-4$, so the weighted residual equations read
\begin{equation}
\int_{-1}^1 R_{n-2}(x) \, (1-x^2) q_{n-4}(x) \, W^{(\gamma)}(x) dx = 0,  \quad \forall \, q_{n-4}(x).
\label{GIG}
\end{equation}
The \emph{Gegenbauer inviscid Galerkin} method is therefore
equivalent to a Gegenbauer tau method with weight $W^{(\gamma+1)}(x)$
and its residual for (\ref{testevp}) is 
\begin{equation}  \label{GIGR}
R_{n-2}(x)=\tau_0 \lambda \, G_{n-2}^{(\gamma+1)}(x) +\tau_1 \lambda \, G_{n-3}^{(\gamma+1)}(x).
\end{equation}
Thus, a Chebyshev (or Legendre) inviscid Galerkin method is a tau method for Chebyshev (or Legendre) polynomials of the \emph{2nd kind}. 
Since we consider a range of the Gegenbauer parameter
$\gamma$, the Gegenbauer tau method also includes some
Gegenbauer Inviscid Galerkin methods.

For completeness, we list the \textit{collocation} approach, where $u_n(x)$ is determined from the
boundary conditions $u_n(\pm 1)=Du_n(\pm 1) = 0$ and enforcing
$R_{n-2}(x_j) = 0$ at the $n-3$ interior Gauss-Lobatto points
$x_j$ such that $DG_{n-2}(x_j) = 0$, $j = 1, . . . , n - 3$,
\cite[\S 2.2]{CHQZ}. The residual (\ref{Rn}) has the form
\cite[eqn.\ (4.5)]{GL83}
\begin{equation}  \label{GCR}
R_{n-2}(x)= (A + B x ) \, DG_{n-2}(x),
\end{equation}
for some $A$ and $B$ independent of $x$. That residual can be written
in several equivalent forms by using the properties of Gegenbauer
polynomials (appendix \ref{GegAppend}). We do not have rigorous results for the collocation method.

\section{Legendre and near-Legendre tau cases}

\label{Legendre}

Here we provide a quicker and more complete derivation of earlier results \cite{DDD98} about spurious eigenvalues for the Legendre and near-Legendre tau case. 
This section provides a useful technical introduction to the problem but is not necessary to derive the main results of this paper. 
Dawkins \etal \cite{DDD98}, focusing only on even modes, use the monomial basis $x^{2k}$  to derive an explicit form for the generalized
eigenvalue problem $A a = \lambda B a$ for the Legendre tau method. In
the monomial basis, the matrix $A$ is upper triangular and
nonsingular, and the matrix $B$ is upper Hessenberg but its first row is identically zero, hence there exists one infinite eigenvalue. 
A perturbation analysis is used to show that the infinite eigenvalue of the
Legendre tau method becomes a large positive eigenvalue 
 for Gegenbauer tau methods with $\gamma < 1/2$ and a large negative
eigenvalue for $ \gamma > 1/2$.

In the Legendre tau method, the polynomial approximation $u_n(x)$ of degree $n$
to problem (\ref{testevp}) satisfies the 4 boundary conditions $u_n=Du_n=0$ at $x=\pm 1$.
Thus $u_n(x)=(1-x^2)^2 p_{n-4}(x)$ and the polynomial $p_{n-4}(x)$ is determined from
 the  weighted residual equations (\ref{jactau}) with $\gamma=1/2$ and $W^{(1/2)}(x)=1$,
\begin{equation}
\int_{-1}^1 \left(\mu D^4 u_n - D^2 u_n \right) \; q_{n-4}(x) dx =
0,  \quad \forall \, q_{n-4}(x).
\label{tauLeg}
\end{equation}
 The mathematical problem is fully specified, except for an
arbitrary multiplicative constant for $u_n(x)$. Choosing various
polynomial bases for $p_{n-4}(x)$ and $q_{n-4}(x)$ will lead to
distinct matrix problems but those problems are all similar to
each other and provide exactly the same eigenvalues, in
exact arithmetic.

 We use the bases $G_l^{(5/2)}(x)$ for $p_{n-4}(x)$ and
$G_k^{(1/2)}(x)$ for $q_{n-4}(x)$, with $k,l=0,\ldots, n-4$, where $\gg_n(x)$ is the Gegenbauer polynomial of degree $n$ for index $\gamma$ (see appendix \ref{GegAppend} and recall that $G_k^{(1/2)}(x)= P_k(x)$ are Legendre polynomials). Thus we write
$u_n(x) = \sum_{l=0}^{n-4} a_l  \; (1-x^2)^2 G_l^{(5/2)}(x), $  
 for some $n-3$ coefficients $a_l$ to be determined.
The tau equations (\ref{tauLeg}) provide the matrix
eigenproblem $\mu A a = B a$, or $\mu \sum_{l=0}^{n-4} A(k,l) a_l
= \sum_{l=0}^{n-4} B(k,l) a_l $ 
 with
\begin{align}
A(k,l) =& \int_{-1}^1 D^4 \left[ (1-x^2)^2 G_l^{(5/2)}(x)\right]
\;
G_k^{(1/2)}(x) \; dx,  \label{Akl}\\
B(k,l)= & \int_{-1}^1 D^2 \left[ (1-x^2)^2 G_l^{(5/2)}(x)\right]
\; G_k^{(1/2)}(x) \;dx, \label{Bkl}
\end{align}
for $k,l=0,\ldots, n-4$. Using (\ref{cl}), these expressions simplify to
\begin{align}
A(k,l) &= \CC_l \int_{-1}^1 \left[D^2 G_{l+2}^{(1/2)}(x) \right] \;
G_k^{(1/2)}(x) \; dx,  \label{Akl2} \\
B(k,l) &= \CC_l \int_{-1}^1 G_{l+2}^{(1/2)}(x) \; G_k^{(1/2)}(x)
\;dx,  \label{Bkl2}
\end{align}
where $\CC_l=\frac{1}{15}(l+1)(l+2)(l+3) (l+4)$.

Since the Legendre polynomials $G_n^{(1/2)}(x)=P_n(x)$ are
orthogonal with respect to the unit weight, equation (\ref{Bkl2}) yields
that $B(k,l) \propto \delta_{k,l+2}$, where $\delta_{k,l+2}$ is the Kronecker delta,  so that $B(0,l)=B(1,l)=0$ for all $l$ and $B$ has non zero elements only on the second subdiagonal.
For $A(k,l)$, use (\ref{dgegrec}) to express $D^2 G_{l+2}^{(1/2)}(x)$ as a linear combination of
$G_{l}^{(1/2)}(x)$, $G_{l-2}^{(1/2)}(x)$, etc.
Orthogonality of the Legendre polynomials $G_n^{(1/2)}(x)$ then implies that
$A(k,l)$ is upper triangular with non-zero diagonal elements.
Hence $A$ is non-singular while the nullspace of $B$ is two-dimensional. The
eigenvalue problem $\mu A a = B a$ therefore has two $\mu=0$
eigenvalues. Since the only non zero elements of $B$ consists of
the sub-diagonal $B(l+2,l)$, the two \emph{right} eigenvectors corresponding
to $\mu=0$ are $a = [0,\ldots,0,1]^T$ and $[0,\ldots,0,1,0]^T$. 
In other words,
\begin{equation}
u_n(x) = (1-x^2)^2 \, G_{n-4}^{(5/2)}(x) \quad \mbox{and} \quad u_n(x) = (1-x^2)^2 \, G_{n-5}^{(5/2)}(x)
\label{spuriousleg}
\end{equation}
satisfy the boundary conditions and the tau equations (\ref{tauLeg}) with $\mu=0=1/\lambda$, for all $n \ge 5$. One mode is even, the other one is odd.
Likewise, the \emph{left} eigenvectors  $b^T =[1,0\ldots,0]$ and
$[0,1,0\ldots,0]$,  satisfy $\mu b^T A = b^T B$ with $\mu =0$. These results are for the Legendre case and $\mu=0$ corresponds to $\lambda = 1/\mu = \infty$.

Now consider the Gegenbauer tau equations for $\gamma-1/2=\epsilon$ with
$|\epsilon| \ll 1$, the near-Legendre case.
The equations are (\ref{tauLeg}) but with the extra weight factor
$W^{(\gamma)}(x) = (1-x^2)^{\epsilon}$ inside the integral. We can figure out what
happens to the $\mu=0$ eigenvalues of the $\epsilon=0$ Legendre case by perturbation. The matrices
$A$ and $B$ and the left and right
eigenvectors, denoted $a$ and $b$ respectively, as well as the eigenvalue $\mu$ now depend on
$\epsilon$. Let
\begin{gather}
A=A_0 + \epsilon A_1 + O(\epsilon^2), \qquad B=B_0 + \epsilon B_1
+ O(\epsilon^2), \\
a=a_0 + \epsilon a_1 + O(\epsilon^2), \qquad b=b_0 + \epsilon b_1
+
O(\epsilon^2), \\
\mu = \mu_0 + \epsilon \mu_1 + O(\epsilon^2),
\end{gather}
where $A_0$ and $B_0$ are the matrices obtained above in (\ref{Akl2}) and (\ref{Bkl2}) for
$\epsilon=0$ while $a_0$ and $b_0^H$ are the corresponding right and left eigenvectors so that
$\mu_0 A_0 a_0 = B_0 a_0$ and $\mu_0 b_0^H A_0 = b_0^H B_0$.
Substituting these $\epsilon$-expansions in the eigenvalue equation
$\mu A a = B a$ and canceling out the zeroth order term, we obtain
\begin{equation}
\mu_1 A_0 a_0 + \mu_0 A_1 a_0 + \mu_0 A_0 a_1= B_1 a_0 + B_0 a_1 +
O(\epsilon).
\end{equation}
Multiplying by $b_0^H$ cancels out the $\mu_0 b_0^H A_0 a_1 =
b_0^H B_0 a_1$ terms so we obtain
\begin{equation}
\mu_1 = \frac{b_0^H B_1 a_0 - \mu_0 b_0^H A_1 a_0}{b_0^H A_0 a_0}.
\label{mu1}
\end{equation}
This expression is general but simplifies further
since we are interested in the perturbation of the zero eigenvalues $\mu_0=0$. This expression for $\mu_1$  is quite simple since $b_0$,
$a_0$ and $A_0$ are the zero$^{th}$ order objects. All we need to
compute when $\mu_0=0$ is the first order correction $B_1$ to the matrix $B$. But  since $a_0$ and $b_0$ have only one non-zero component as given at the end of the previous paragraph, we only need to calculate two components of the $B$ matrix.
For $n$ even, all that is needed are the first order corrections to $B(0,n-4) $ for the even mode and to $B(1,n-5)$ for the odd mode. For $n$ odd,  we need $B(0,n-5)$ for the even mode and $B(1,n-4)$ for the odd mode, however since even and odd modes decouple in this problem, it suffices to compute both even and odd modes in only one case of $n$ even or odd.
The matrix elements in the $\epsilon \ne 0$ cases are still given by (\ref{Akl2}) and (\ref{Bkl2}) but with the extra $(1-x^2)^{\epsilon}$ weight factor inside the integrals. Since
  $G_{1}^{(1/2)}(x)=x$ and $G_{n}^{(1/2)}(x)=P_n(x)$, the Legendre
polynomial of degree $n$, we obtain
\begin{align}
B(0,n-4) = & \CC_{n-4} \int_{-1}^1P_{n-2}(x)
\;(1-x^2)^{\epsilon} \;dx = \epsilon B_1(0,n-4) + O(\epsilon^2), \label{B1n-4}\\
B(1,n-5) = & \CC_{n-5}\int_{-1}^1 x \; P_{n-3}(x)
\; (1-x^2)^{\epsilon} \;dx = \epsilon B_1(1,n-5) + O(\epsilon^2), \label{B1n-5}
\end{align}
with $\CC_l$ as defined in (\ref{Bkl2}).
The integrals are readily evaluated and details are provided in appendix  \ref{Asymptotics}. Using
(\ref{B1}), (\ref{Kn2}) and (\ref{An-4}) we obtain for the
even mode (for $n$ even) that
\begin{equation}
\mu_1=\frac{B_1(0,n-4)}{A_0(0,n-4)} = \frac{-4}{(n-2)^2(n-1)^2}.
\label{mu1even}
\end{equation}
This matches the formula in Dawkins \etal \cite[page 456]{DDD98} since their
$2N=n-4$ and $2 \nu-1=2 \epsilon$. Likewise, using
 (\ref{B2eval}), and (\ref{An-5}) for the odd mode (with $n$ even) yields
\begin{equation}
\mu_1=\frac{B_1(1,n-5)}{A_0(1,n-5)} = \frac{-4}{(n-4)^2(n-1)^2}.
\label{mu1odd}
\end{equation}
Again, if $n$ is odd then $\mu_1$ for the even mode is given by
(\ref{mu1even}) but with $n-1$ in lieu of $n$. Likewise for $n$ odd,
the odd mode is given by (\ref{mu1odd}) with $n+1$ in lieu of $n$.
Finally, since $\lambda= 1/\mu$, the $\lambda = \infty$ eigenvalues
in the Legendre tau case, become $\lambda = 1/(\epsilon \mu_1 +
O(\epsilon^2)) \sim 1/(\epsilon \mu_1)$ in the near-Legendre cases.
From (\ref{mu1even}) and (\ref{mu1odd}), these eigenvalues will be
$O(n^4/\epsilon)$. Furthermore they will be positive when $\epsilon
< 0$ (\ie \emph{spurious} when $\gamma < 1/2$) but negative when
$\epsilon
>0$.

\section{Characteristic Polynomials}

\label{polys}

 For the model problem (\ref{testevp}), we can bypass
the matrix eigenproblem of section \ref{Legendre} to directly derive
the characteristic
polynomial for the eigenvalues $\mu=1/\lambda$. To do so, invert equation
(\ref{Rn}) to express the polynomial approximation $D^2u_n(x)$ in
terms of the residual $R_{n-2}(x)$
\begin{equation}  \label{D2u_n}
D^2u_n(x)=\mu \sum_{k=0}^{\infty} \mu^k D^{2k} R_{n-2}(x)
\end{equation}
where $\mu = 1/\lambda$.
 The inversion (\ref{D2u_n}) follows from
application of the geometric (Neumann) series for
$(1-\mu D^2)^{-1}= \sum_{k=0}^{\infty} \mu^k D^{2k}$ which terminates
since $R_{n-2}(x)$ is a polynomial.
Thus, $u_n(x)$ can be
computed in terms of the unknown tau coefficients by double integration of (\ref{D2u_n}) and application of the boundary conditions. We can assume that $\lambda \ne 0$ because $\lambda=0$
with $u_n(\pm 1)=Du_n(\pm 1)=0$ necessarily corresponds to the
trivial solution $u_n(x)=0$, $\forall x$ in $[-1,1]$, as noted in the previous section.

The Gegenbauer polynomials are even in $x$ for $n$ even and odd for $n$ odd
 (\ref{Goddeven}).
 The symmetry of the differential equation (\ref{testevp}) and of the Gegenbauer
polynomials allows decoupling of the discrete problem into even
and odd solutions. This parity reduction leads to simpler
residuals and simpler forms for the corresponding characteristic
polynomials. The residual in the parity-separated Gegenbauer case
contains only one term
\begin{equation}
R_{n-2}(x)=\tau_0 \lambda \, G_{n-2}^{(\gamma)}(x),
\label{RnGeg}
\end{equation}
instead of (\ref{GTR}), where $G_n^{(\gamma)}(x)$ is the Gegenbauer polynomial of degree $n$
with $n$  even for even solutions and odd for odd solutions.
Substituting (\ref{RnGeg}) in (\ref{D2u_n}) and renormalizing $u_n(x)$ by $\tau_0$ gives
\begin{equation}
D^2u_n(x)=\sum_{k=0}^{\infty} \mu^k D^{2k}G^{(\gamma)}_{n-2}(x).
\label{D2u}
\end{equation}
For $\gamma > 1/2$,  the identity (\ref{gege2}) in the form
$2 \gamma G_{n-2}^{(\gamma)}(x) =DG_{n-1}^{(\gamma-1)}(x)$  can be used to write
(\ref{D2u}) in the form
\begin{equation}
D^2u_n(x)=\frac{1}{2\gamma} \sum_{k=0}^{\infty} \mu^k D^{2k+1}G^{(\gamma-1)}_{n-1}(x),
\label{D2utoo}
\end{equation}
which integrates to
\begin{equation}
Du_n(x)=\frac{1}{2\gamma} \sum_{k=0}^{\infty} \mu^k D^{2k}G^{(\gamma-1)}_{n-1}(x) + C,
\label{Du}
\end{equation}
where $C$ is an arbitrary constant.

\subsection{Even Solutions}

For even solutions $u_n(x) = u_n(-x)$,  $n$ is even and $Du_n(x)$ is odd so $C=0$ in (\ref{Du}). The boundary condition $Du_n(1)=0$ gives the characteristic equation for $\mu$ (for $n$ even and $\gamma>1/2$)
\begin{equation}
\label{GTCPE2}
\sum_{k=0}^{\infty} \mu^k \, D^{2k}G^{(\gamma-1)}_{n-1}(1)=0.
\end{equation}

\subsection{Odd Solutions}

For odd solutions,  $u_n(x)=-u_n(-x)$, $n$ is odd and the boundary condition $Du_n(1)=0$ requires that
\begin{equation}
C=- \frac{1}{2 \gamma} \sum_{k=0}^{\infty}\mu
^{k}D^{2k}G_{n-1}^{(\gamma-1 )}(1).
\label{Codd}
\end{equation}
Substituting this $C$ value in (\ref{Du}) and integrating gives
\begin{equation}
2\gamma \,u_{n}(x)=\sum_{k=0}^{\infty}\mu ^{k}D^{2k-1}G_{n-1}^{(\gamma-1)}(x)-
x \sum_{k=0}^{\infty}\mu^{k}D^{2k}G_{n-1}^{(\gamma-1)}(1)  \label{equo}
\end{equation}
where we must define
\begin{equation}
D^{-1}G_{n-1}^{(\gamma-1 )}(x)=\int_{0}^{x}G_{n-1}^{(\gamma-1
)}(s)\;ds =  \frac{G_n^{(\gamma-1)}(x)-G_{n-2}^{(\gamma-1)}(x)}{2(n+\gamma-2)}
\label{Dm1odd}
\end{equation}
since $u_n(x)$ and $n$ are odd, where we have used (\ref{dgegrec}) to evaluate the integral
and the symmetry (\ref{Goddeven}) so that $G_n(0)=G_{n-2}(0)=0$ for $n$ odd.
The boundary condition $u_n(1)=0$ yields the
characteristic polynomial equation (for $n$ odd and $\gamma>1/2$)
\begin{equation}
\sum_{k=0}^{\infty} \mu^k D^{2k-1}G^{(\gamma-1)}_{n-1}(1)-
\sum_{k=0}^{\infty} \mu^k D^{2k}G^{(\gamma-1)}_{n-1}(1)=0.
\label{GTCPO}
\end{equation}

For $\gamma  > 3/2$, we can use identity (\ref{gege2}) in the form
$2 (\gamma-1) G_{n-1}^{(\gamma-1)}(x) =DG_{n}^{(\gamma-2)}(x),$ to write the characteristic equation
(\ref{GTCPO}) as
\begin{equation}
\sum_{k=0}^{\infty} \mu^k D^{2k}G^{(\gamma-2)}_{n}(1)-
\sum_{k=0}^{\infty} \mu^k D^{2k+1}G^{(\gamma-2)}_{n}(1)=0.
\label{GTCPO2}
\end{equation}
For $1/2 < \gamma \le 3/2$, this cannot be used since $\gamma-2 < -1/2$, but using (\ref{Dm1odd}) for the $D^{-1}$ term in the first sum, the characteristic equation (\ref{GTCPO}) can be written
\begin{equation}
\mu \sum_{k=0}^{\infty} \mu^k D^{2k+1}G^{(\gamma-1)}_{n-1}(1)
- \frac{G_{n-2}^{(\gamma-1)}(1)-G_{n}^{(\gamma-1)}(1)}{2(n+\gamma-2)}
-\sum_{k=0}^{\infty} \mu^k D^{2k}G^{(\gamma-1)}_{n-1}(1)=0.
\label{GTCPO3}
\end{equation}

\section{Zeros of Characteristic Polynomials}

\label{Zeros}

Here we prove that the zeros of the characteristic polynomial equations (\ref{GTCPE2}) and (\ref{GTCPO}) are real, negative, and distinct for $1/2 < \gamma \le 7/2$. Some background material is needed.

\subsection{Stable Polynomials and the Hermite Biehler Theorem}

A polynomial $p(z)$ is \textit{stable} if and only if all its zeros have negative real parts.
Stable polynomials can arise as characteristic polynomials of a
numerical method applied to a differential equation as in \cite{CCW} for $Du=\lambda u$ with $u(1)=0$ and in other dynamical systems applications.
The characterization of stable polynomials that is most useful here is given by \cite{holtz}, \cite[p.197]{RS},
\begin{theorem}
{The Hermite-Biehler Theorem.} \label{HBthm}  The real polynomial
$p(z)\,=\,\Omega(z^2)+z \Theta(z^2)$
is stable if and only if $\Omega(\mu)$ and $\Theta(\mu)$ form a
positive pair.
\end{theorem}

\begin{definition} Two real polynomials $\Omega(\mu)$ and $\Theta(\mu)$  of
degree $n$ and $n-1$ (or $n$) respectively, form a positive pair if:  

(a) the roots $\mu_1,\mu_2,\ldots,\mu_{_{n}} $ of $\Omega(\mu)$ and
$\mu_1^{^{\prime}},\mu_2^{^{\prime}},\ldots,\mu_{_{n-1}} ^{^{\prime}}$
 (or $ \mu_1^{^{\prime}},\mu_2^{^{\prime}},\ldots,\mu_{_{n}} ^{^{\prime}}$)
of $\Theta(\mu)$ are all distinct, real and negative.

(b) the roots interlace as follows: 
$\mu_1<\mu_1^{^{\prime}}<\mu_2<\cdots <\mu_{_{n-1}} ^{^{\prime}}<\mu_{_{n}} <0$ ( or $\mu_1^{^{\prime}}<\mu_1<\cdots<\mu_{_{n}} ^{^{\prime}}<\mu_{_{n}} <0$ )

(c) the highest coefficients of $\Omega(\mu)$ and $\Theta(\mu)$
are of like sign.
\end{definition}

We will use the following theorem about positive pairs
\cite[p198]{L},

\begin{theorem}
\label{thm3}
Any nontrivial real linear combination of two
polynomials of degree $n$ (or $n$ and $n-1$)  with interlacing roots has real roots.
\end{theorem}
(Since such a linear combination changes sign $n-1$ times along the real axis, it has $n-1$ real roots.
Since it is a real polynomial of degree $n$,  the remaining root is real also.)

\subsection{Eigenvalues for even modes}
\label{zeroseven}

In \cite{CCW} and \cite{CW} we study the Gegenbauer tau method for
$D^2u=\lambda u$ with $u(\pm 1)=0$, which leads to the
characteristic polynomials
$\sum_{k=0}^{\infty} \mu^k D^{2k} G_n^{(\gamma)}(1)$.
The derivation of that result is entirely similar to that in sections
\ref{residual} and \ref{polys}. The strategy to prove that the Gegenbauer tau method for that 2nd order problem has real,  negative and distinct roots is to show stability of the polynomial
\begin{equation}  \label{stable}
p(z) =\sum_{k=0}^{n} z^k \, D^{k} G_n^{(\gamma)}(1)
\end{equation}
for $-1/2 < \gamma \le 3/2$ then to use the Hermite Biehler Theorem  to deduce that
\begin{theorem}
\label{pp}
For  $-1/2 < \gamma \le {3}/{2}$, the polynomials
\begin{equation}
\Omega_n^{(\gamma)}(\mu)=\sum_{k=0}^{\infty} \mu^k D^{2k}
G_n^{(\gamma)}(1) \qquad \mbox{and} \qquad
\Theta_n^{(\gamma)}(\mu)=\sum_{k=0}^{\infty} \mu^k D^{2k+1}
G_n^{(\gamma)}(1)
\label{OmThet2}
\end{equation}
form a positive pair.  From (\ref{gege2}) this is equivalent to stating that the
polynomials
\begin{equation}
\Omega_n^{(\gamma)}(\mu)=\sum_{k=0}^{\infty} \mu^k D^{2k}
G_n^{(\gamma)}(1) \qquad \mbox{and} \qquad
\Omega_{n-1}^{(\gamma+1)}(\mu)=\sum_{k=0}^{\infty} \mu^k D^{2k}
G_{n-1}^{(\gamma+1)}(1)
\label{OmThet2too}
\end{equation}
also form a positive pair. Combining the $\gamma$ and $\gamma+1$ ranges in (\ref{OmThet2too}) yields that $\Omega_n^{(\gamma)}(\mu)$ has real, negative, and distinct roots for
 $-1/2 < \gamma \le 5/2$.
\end{theorem}
Stability of (\ref{stable}) is proven in \cite[Theorem 1]{CCW} for the broader class of Jacobi
polynomials. The basic ideas of the proof are along the lines of Gottlieb \cite{G81} and Gottlieb and Lustman's work \cite{GL83} on stability of the Chebyshev collocation method for the 1st and 2nd order operator.

For the 4th order problem (\ref{testevp}), $D^4 u = \lambda D^2u$ with $u(\pm 1)= Du(\pm 1)=0$,
the Gegenbauer tau method gives the characteristic polynomial (\ref{GTCPE2})
for even solutions. This is the polynomial $\Omega_{n-1}^{(\gamma-1)}(\mu)$  of (\ref{OmThet2too}) that appears for the 2nd order problem \cite{CW,CCW}  and is known to have real, negative, and distinct eigenvalues for $-1/2 < \gamma -1 \le 5/2$, that is for $1/2 < \gamma \le 7/2$.
Hence, it follows directly from (\ref{GTCPE2}) and theorem \ref{pp} that
 the Gegenbauer tau approximation for even solutions of
problem (\ref{testevp}), $D^4 u = \lambda D^2 u$ with $u(\pm 1) = Du(\pm 1)=0$, has real, negative, and distinct eigenvalues for $1/2 < \gamma \le 7/2$.

This result is sharp.  For $\gamma>7/2$ and sufficiently large $n$, our numerical computations show that the polynomial has a pair of complex eigenvalues. For $\gamma < 1/2$ the polynomial (\ref{GTCPE2}) has a real positive eigenvalue as first proven in \cite{DDD98}. The proof  goes as follows. For $\gamma < 1/2$,  we cannot use (\ref{Du}) since  $\gamma -1$ is below the range of definition of Gegenbauer polynomials (appendix \ref{GegAppend}).
Instead, integrate (\ref{D2u}) and use identity (\ref{dgegrec}) to obtain
$\int_0^1 G_{n-2}^{(\gamma)} dx =
(G_{n-1}^{(\gamma)}(1)- G_{n-3}^{(\gamma)}(1))/(2(\gamma+n-2))$ since
$G_{n-1}(0)=G_{n-3}(0)=0$ for $n$ even (\ref{Goddeven}).
Using formula (\ref{Gat1}), one shows that
$G_{n}^{(\gamma)}(1)$ increases with $n$ if $\gamma>1/2$ but
decreases with $n$ if $\gamma<1/2$. Hence, the constant term is
negative for $\gamma<1/2$ while all the other coefficients of the
characteristic polynomial can be shown to be positive using (\ref{gege2}) and (\ref{Gat1}). Therefore, there is one
real positive eigenvalue as proved in \cite{DDD98}. For $\gamma=1/2$, the Legendre case, the constant term vanishes and there is a $\mu=0$ ($\lambda=\infty$) eigenvalue as established in section \ref{Legendre}.

\begin{rmk}
Exact even solutions of  $D^4 u = \lambda D^2 u$ with $u(\pm 1)=Du(\pm 1)=0$ obey $D^3 u = \lambda Du + C$ with $C=0$, since $D^3u$ and $Du$ are odd. Thus even solution eigenvalues of (\ref{testevp}) are equal to the eigenvalues for  odd solutions of the 2nd order problem $D^2 w = \lambda w$ with $w=0$ at $x=\pm 1$, with $w=Du$. The Gegenbauer tau version of this property is that eigenvalues for even Gegenbauer tau solutions of (\ref{testevp}) of order $n$ (even) and index $\gamma$  are equal to the eigenvalues for odd Gegenbauer tau solutions of $D^2 w = \lambda w$, $w(\pm 1)=0$ of order $n-1$ (odd) and index $\gamma-1$. This follows directly from (\ref{GTCPE2}) and \cite{CCW,CW}.
\end{rmk}

 \subsection{Eigenvalues for odd modes}

 \label{zerosodd}

The reduction of the 4th order problem (\ref{testevp}) to the 2nd order problem  does not hold for odd modes which have the characteristic equation (\ref{GTCPO}). In fact, all previous theoretical work focused only on the even modes \cite{DDD98,GO77}. The same general strategy as for the even case  led us to prove the stability of a shifted version of polynomial (\ref{stable}).

\begin{theorem}
\label{thm41} Let $G^{(\gamma)}_n(x)$ denote the non standard
Gegenbauer polynomial of degree $n$ as defined in appendix \ref{GegAppend}, then the polynomial
\begin{equation}
p^{(\gamma)}_n(z)=
 \frac{G_{n-1}^{(\gamma)}(1) - G_{n+1}^{(\gamma)}(1)}{2(n+\gamma)}
+ \sum_{k=0}^n z^k \, D^k G_n^{(\gamma)}(1)
\end{equation}
is stable for $-1/2 < \gamma \le 1/2$.
\end{theorem}

The proof is elementary but technical, it is given in Appendix \ref{ProofAppend}.
We also need the following simple lemma. This lemma will help us
determine the sign of the coefficients of the characteristic
polynomials.

\begin{lemma}
\label{lem1} With $G_n^{(\gamma)}(x)$ as defined in appendix \ref{GegAppend},
 the expression
\begin{equation}
D^{k+1}G_n^{(\gamma)}(1)-D^{k}G_n^{(\gamma)}(1) \ge 0
\end{equation}
 for $k=0, \dots, n-1$ and $\gamma > -1/2$.
\end{lemma}

\begin{proof}
From (\ref{dkG1})
\begin{equation}
D^{k}G_n^{(\gamma)}(1)=\frac{2^{k-1} \Gamma{(\gamma +k)}
\Gamma{(n+2\gamma +k)}}{(n-k)!\Gamma{(\gamma+1)} \Gamma{(2\gamma
+2k)}},
\end{equation}
where $\Gamma(z)$ is the standard gamma function.
For $k>0$ and given that $\gamma > -1/2$, the sign of the above
expression is positive since all individual terms are positive. In
the $k=0$ case, the
sign of the expression is determined by the terms
 $\frac{\Gamma{(\gamma)}}{\Gamma{(2\gamma)}}$
 since all other terms are positive. If $-1/2 <
\gamma <0$ both numerator and denominator are negative and thus
their ratio is positive. If $\gamma >0 $ the two terms are
positive and thus again their ratio is positive. For $\gamma=0$,
 a simple limiting argument shows positiveness again, in fact from (\ref{gegchebleg}),  $G_n^{(0)}(1) = T_n(1)/n = 1/n$.

Taking the ratio $D^{k+1}G_n^{(\gamma)}(1)/D^{k}G_n^{(\gamma)}(1)$
and making some simplifications gives
\begin{equation}
\frac{D^{k+1}G_n^{(\gamma)}(1)}{D^{k}G_n^{(\gamma)}(1)}\,=
\,\frac{(2\gamma+n+k)(n-k)}{(2\gamma+2k+1)}.
\end{equation}
Since $k \le n-1 $ then $2\gamma+2k+1 \le 2\gamma+(n-1)+k+1 =
2\gamma+n+k $. Thus
\begin{equation}
\frac{D^{k+1}G_n^{(\gamma)}(1)}{D^{k}G_n^{(\gamma)}(1)}\,\ge 1
\quad k=0 \dots n-1
\end{equation}
and since both derivatives are positive, the lemma follows.
\end{proof}

We now have all the tools to prove

\begin{theorem}
\label{thmtau} The Gegenbauer tau approximation to problem
(\ref{testevp}) has  real, negative, and distinct eigenvalues for
$1/2 < \gamma \le 7/2$. This $\gamma$ range is sharp, spurious positive eigenvalues exist for $\gamma<1/2$ and complex eigenvalues arise for $7/2< \gamma$.
\end{theorem}

\begin{proof}
This has already been proven in section \ref{zeroseven} for even solutions.
For odd solutions, we need to consider two separate cases.

\textbf{Case 1.} {$3/2 < \gamma \le 7/2$}. The characteristic polynomial
(\ref{GTCPO2})
\begin{equation}
\sum_{k=0}^{\infty}\,\mu^k\,D^{2k}\,G_{n}^{(\gamma-2)}(1)\,-\,
\sum_{k=0}^{\infty}\,\mu^k\,D^{2k+1}\,G_{n}^{(\gamma-2)}(1)\,=
\Omega^{(\gamma-2)}_n(\mu)-\Theta^{(\gamma-2)}_{n}(\mu)
\end{equation}
is a linear combination of the polynomials $\Omega_n^{(\gamma-2)}(\mu)$ and
$\Theta_n^{(\gamma-2)}(\mu)$, which form a positive pair for
$3/2<\gamma \le 7/2$, by theorem \ref{pp}. Therefore, by theorem
\ref{thm3} this characteristic polynomial has real roots. Then by
lemma \ref{lem1} we deduce that all its coefficients are of the same sign,
 hence all its roots must be negative.

\textbf{Case 2.} {$1/2 < \gamma \le 3/2$}.
The polynomial
\begin{equation}
p^{(\gamma-1)}_{n-1}(z) = \Lambda(z^2) + z \Phi(z^2)
\end{equation}
with $p^{(\gamma)}_{n}(z)$ as in theorem \ref{thm41}, is stable for the
desired range of parameters by theorem  \ref{thm41}, so
 the Hermite Biehler theorem (theorem \ref{HBthm})  implies that the polynomials
\begin{multline}
\Lambda(\mu)=
\frac{G_{n-2}^{(\gamma-1)}(1) - G_{n}^{(\gamma-1)}(1)}{2(n+\gamma-2)}
+  \sum_{k=0}^{\infty} \mu^k \, D^{2k} G_{n-1}^{(\gamma-1)}(1)
 \\
=
\frac{G_{n-2}^{(\gamma-1)}(1) -G_{n}^{(\gamma-1)}(1)}{2(n+\gamma-2)}
+ \Omega_{n-1}^{(\gamma-1)}(\mu),
\end{multline}
and
\begin{equation}
\Phi(\mu)=\sum_{k=0}^{\infty} \mu^k \, D^{2k+1} G_{n-1}^{(\gamma-1)}(1)
=\Theta_{n-1}^{(\gamma-1)}(\mu)
\end{equation}
form a positive pair, with $\Omega_n^{(\gamma)}(\mu)$ and
$\Theta_n^{(\gamma)}(\mu)$ as defined in theorem \ref{pp}. Thus $\mu
\Phi(\mu)$ and  $\Lambda(\mu)$ have interlacing roots and any real
linear combination of the two must have real roots (theorem
\ref{thm3}). Now the characteristic polynomial (\ref{GTCPO3}) is in
fact the linear combination $\mu \Phi(\mu) - \Lambda(\mu)$ so it has real roots.
Its constant term is equal to
\begin{equation}
 \frac{G_{n}^{(\gamma-1)}(1) -G_{n-2}^{(\gamma-1)}(1)}  {2(n+\gamma-2)}
 -G_{n-1}^{(\gamma-1)}(1)
\end{equation}
which is negative if $-1/2 < \gamma -1 \le 1/2 $, that is $1/2 < \gamma \le 3/2$,  from (\ref{Gat1}).
 All other coefficients of
the characteristic polynomial $\mu \Phi(\mu) - \Lambda(\mu)$  are negative by lemma \ref{lem1}.
Since all coefficients have the same sign and all roots are real, all the roots must be negative.
\end{proof}

\section{Galerkin and Modified tau methods}
\label{numerical}

Our main theorem \ref{thmtau} can be expressed in terms of the inviscid Galerkin and Galerkin methods since these methods are equivalent to Gegenbauer tau methods with index $\gamma+1$ and $\gamma+2$, respectively, as shown in section \ref{residual}.

\begin{corollary}
\label{GIGthm} The Gegenbauer Inviscid Galerkin approximation
 to problem (\ref{testevp})  has real, negative, and distinct eigenvalues for $-1/2 < \gamma \le 5/2$.
\end{corollary}

\begin{corollary}
\label{GGthm} The Gegenbauer Galerkin approximation 
to problem (\ref{testevp}) has real negative eigenvalues for $-1/2< \gamma \le 3/2$.
\end{corollary}

\begin{corollary}
Since Chebyshev corresponds to $\gamma=0$ and Legendre to $\gamma=1/2$, the Chebyshev and Legendre tau approximations to problem (\ref{testevp})  have spurious eigenvalues, but the Chebyshev and Legendre inviscid Galerkin ($\gamma=1$ and $3/2$, respectively) and Galerkin $(\gamma=2$ and $5/2$, respectively) approximations provide real, negative, and distinct eigenvalues.
\end{corollary}

Finally, we prove that the modified tau method introduced by Gottlieb and Orszag \cite{GO77} and developed by various authors \cite{GTD89}, \cite{MMB90} is equivalent to the Galerkin method. McFadden, Murray and Boisvert \cite{MMB90} have already shown the equivalence between the modified Chebyshev tau and the Chebyshev Galerkin methods. Our simpler proof generalizes their results to the Gegenbauer class of approximations.

The idea for the modified tau method, widely used for time-marching, starts with the substitution
$v(x)=D^2u(x)$. Problem (\ref{testevp}) reads
\begin{equation}
 D^2u=v, \qquad D^2v = \lambda v,  \quad \mbox{with} \; u=Du=0 \;\; \mbox{at} \;
x=\pm 1. \label{modeq}
\end{equation}
If we approximate $u(x)$ by a polynomial of degree $n$, then $v=D^2 u$ suggests that $v$ should be a polynomial of degree $n-2$, however the modified tau method approximates both $u(x)$ and $v(x)$ by polynomials of degree $n$,
\begin{equation}
\begin{split}  \label{equv}
u_{n}(x)=\sum_{k=0}^n \hat{u}_k G_k^{(\gamma)}(x) , \qquad
D^2u_{n}(x)=\sum_{k=0}^{n-2} \hat{u}_k^{(2)} G_k^{(\gamma)}(x), \\
v_{n}(x)=\sum_{k=0}^n \hat{v}_k G_k^{(\gamma)}(x) , \qquad
D^2v_{n}(x)=\sum_{k=0}^{n-2} \hat{v}_k^{(2)} G_k^{(\gamma)}(x), \\
\end{split}
\end{equation}
where the superscripts indicate the Gegenbauer coefficients of the
corresponding derivatives. These can be expressed in terms of the Gegenbauer coefficients of the original function using (\ref{dgegrec}) twice, as in the Chebyshev tau method \cite{CHQZ},  \cite{O71}. Hence there are $2n+2$ coefficients to be determined, $\hat{u}_0,\ldots,\hat{u}_n$,
 $\hat{v}_0,\ldots,\hat{v}_n$. In the modified tau method, these are determined by the 4 boundary conditions and the $2n-2$ tau equations obtained from orthogonalizing the residuals of both equations $D^2 u = v$ and $D^2 v = \lambda v$ to the first $n-1$ Gegenbauer polynomials
$G_0^{(\gamma)}(x), \ldots, G_{n-2}^{(\gamma)}(x)$ with respect to the Gegenbauer weight $(1-x^2)^{\gamma-1/2}$. In terms, of the expansions (\ref{equv}), these weighted residual equations have the simple form
\begin{equation}
 \label{modeq2}
 \begin{aligned}
\hat{u}_k^{(2)}=& \hat{v}_k, & \quad 0 \le k \le n-2, \\
\hat{v}_k^{(2)}=& \lambda \hat{v}_k, & \quad 0 \le k \le n-2.
\end{aligned}
\end{equation}
McFadden \etal \cite{MMB90} showed that the modified Chebyshev tau
 is equivalent to the Chebyshev Galerkin method for this
particular problem.  We provide a simpler proof for the general setting of
 Gegenbauer polynomials.

\begin{theorem}
The modified Gegenbauer tau method proposed in \cite{GO77}  is equivalent to the Gegenbauer Galerkin method for problem (\ref{testevp}).
\end{theorem}

\begin{proof}
Let the polynomial approximations and their derivatives as in
(\ref{equv}). Then the tau equations (\ref{modeq2}) are equivalent to the residual equations
\begin{equation}
\begin{aligned}
v_n(x)-D^2u_n (x)= &\; \hat{v}_{n-1} G_{n-1}^{(\gamma)}(x)
 &+& \; \hat{v}_n G_{n}^{(\gamma)}(x), \\
 (\lambda - D^2)v_n(x) = & \lambda \hat{v}_{n-1} G_{n-1}^{(\gamma)}(x)
&+& \lambda \hat{v}_n G_{n}^{(\gamma)}(x).
\end{aligned}
\end{equation}
Combining the two yields
\begin{equation}
(\lambda - D^2)D^2u_n (x)=\hat{v}_{n-1} \,D^2G_{n-1}^{(\gamma)}(x)
+ \hat{v}_n \,D^2G_{n}^{(\gamma)}(x),
\quad u_n (\pm 1)=Du_n (\pm 1)=0
\end{equation}
which using (\ref{gege2}) is equivalent to
\begin{equation}
(\lambda - D^2)D^2u_n (x)=\tau_0 \lambda G_{n-2}^{(\gamma+2)}(x)
+\tau_1 \lambda G_{n-3}^{(\gamma+2)}(x), \quad u_n (\pm 1)=Du_n
(\pm 1)=0.
\end{equation}
with $\hat{v}_n= 4 \lambda (\gamma+1)(\gamma+2) \tau_0$ and
$\hat{v}_{n-1} =4 \lambda(\gamma+1)(\gamma+2) \tau_1$. This
is exactly the Gegenbauer Galerkin method, as
given in equation (\ref{GGR}). 
\end{proof}

Since the modified Gegenbauer tau is equivalent to the Galerkin method,
 the results in section \ref{Zeros} imply that the
modified tau method for problem \ref{testevp} has real and
negative eigenvalues for $-1/2 < \gamma \le 3/2$. This includes Chebyshev for $\gamma=0$ and Legendre for $\gamma=1/2$.

\section*{Acknowledgments}
We thank Jue Wang for helpful calculations in the early stages of
this work. This research was partially supported by NSF grant
DMS-0204636.

\appendix  

\section{Gegenbauer (Ultraspherical) Polynomials}

\label{GegAppend}

The Gegenbauer (a.k.a.\ Ultraspherical) polynomials $C_n^{(\gamma)}(x)$,
$\gamma > -1/2$, of degree $n$ are the Jacobi polynomials with
 $\alpha=\beta=\gamma-1/2$, up to normalization \cite[22.5.20]{AS}.
They are symmetric (even for $n$ even and odd for $n$ odd)
orthogonal polynomials with weight function
$W^{(\gamma)}(x)=(1-x^2)^{\gamma-1/2}$. Since the standard
normalization \cite[22.3.4]{AS}, is singular for the Chebyshev case
$\gamma=0 $, we use a non-standard normalization that includes the
Chebyshev case but preserves the simplicity of the Gegenbauer
recurrences. Set
\begin{equation}
G_0^{(\gamma)}(x):=1, \qquad G_n^{(\gamma)}(x) :=
 \frac{C_n^{(\gamma)}(x)} {2\gamma}, \quad n \ge 1.
\end{equation}
These Gegenbauer polynomials satisfy the orthogonality relationship
\begin{equation}
\int_{-1}^{1} (1-x^2)^{\gamma-1/2} \,G_m^{(\gamma)} G_n^{(\gamma)}
\,dx
=\left\{%
\begin{array}{ll}
0, & m \ne n, \\
h_n^{\gamma}, & m = n,
\end{array}%
\right.  \label{gegortho}
\end{equation}
where \cite[22.2.3]{AS},
\begin{equation}
h_0^{\gamma}=\frac{\pi 2^{-2 \gamma} \Gamma(2 \gamma+1)}{
\Gamma^2(\gamma+1)} \qquad \mbox{and} \qquad h_n^{\gamma} =
\frac{\pi 2^{-1-2 \gamma} \Gamma(n+2 \gamma)}{(n+\gamma)n!
\Gamma^2(\gamma+1)} \qquad n \ge 1.
\end{equation}
The Sturm-Liouville form of the Gegenbauer equation for
$G^{(\gamma)}_n(x)$ is
\begin{equation}
D \left[(1-x^2)^{\gamma+1/2}\, D G^{(\gamma)}_n(x) \right]=-
n(n+2\gamma) \, (1-x^2)^{\gamma-1/2} \,G^{(\gamma)}_n(x)
\label{gegsl}
\end{equation}
They satisfy the derivative recurrence formula
\begin{equation}  \label{gege2}
\frac{d}{dx} G_{n+1}^{(\gamma)} = 2 (\gamma+1) G_n^{(\gamma+1)},
\end{equation}
(for $C_n^{(\gamma)}$ this is formula \cite[A.57]{Boyd}), and
their three-term recurrence takes the simple form
\begin{equation}
(n+1)G_{n+1}^{(\gamma)}= 2(n+\gamma) x G_{n}^{(\gamma)} - (n-1 + 2
\gamma) G_{n-1}^{(\gamma)}, \quad n\ge 2, ,  \label{gegrec}
\end{equation}
with 
\begin{equation}
G_{0}^{(\gamma)}(x) =1, \quad G_{1}^{(\gamma)}(x) =x, \quad
G_2^{(\gamma)}=(\gamma+1)x^2-\frac{1}{2}.
\end{equation}
These recurrences can be used to verify the odd-even symmetry of Gegenbauer polynomials
\cite[22.4.2]{AS}
\begin{equation}
G_n^{(\gamma)}(x) = (-1)^n G_n^{(\gamma)}(-x).
\label{Goddeven}
\end{equation}
Differentiating the recurrence (\ref{gegrec}) with respect to $x$
and
subtracting from the corresponding recurrence for $\gamma+1$ using (\ref%
{gege2}), yields \cite[22.7.23]{AS}
\begin{equation}  \label{ggp}
(n+\gamma) G_{n}^{(\gamma)} = (\gamma+1) \left[G_{n}^{(\gamma+1)}
- G_{n-2}^{(\gamma+1)} \right], \qquad n \ge 3.
\end{equation}
Combined with (\ref{gege2}), this leads to the important
derivative recurrence between Gegenbauer polynomials of same
index $\gamma$
\begin{equation}
\begin{split}
G_0^{(\gamma)}(x) = D G_1^{(\gamma)}(x),& \quad 2 (1+\gamma)
G_1^{(\gamma)}(x) = D G_2^{(\gamma)}(x), \\
2 (n+\gamma) G_{n}^{(\gamma)} =\; & \frac{d}{dx}
\left[G_{n+1}^{(\gamma)} - G_{n-1}^{(\gamma)} \right].
\label{dgegrec}
\end{split}
\end{equation}
Evaluating the Gegenbauer polynomial at $x=1$ we find
\cite[22.4.2]{AS},
\begin{equation}  \label{Gat1}
\begin{aligned}
G_n^{(\gamma)}(1)=\frac{1}{2\gamma}C_n^{(\gamma)}(1)=\frac{1}{2\gamma}
\binom{2\gamma + n-1}{n}
=& \frac{(2\gamma+n-1)(2\gamma+n-2)\cdots (2\gamma+1)}{n!} \\
=&  \frac{\Gamma{(2\gamma + n)}}{n! \, \Gamma{(2\gamma+1)}}
\end{aligned}
\end{equation}
for $n\ge 2$, with $G_1^{(\gamma)}(1)=G_0^{(\gamma)}(1)=1$, where $\Gamma(z)$ is the standard gamma function \cite{AS}. Note that $G_n^{(\gamma)}(1)>0 $ for $\gamma>-1/2$ and that it decreases with increasing $n$ if  $-1/2 < \gamma < 1/2$ but increases with $n$ if $1/2 < \gamma$.

Now (\ref{gege2}) gives
\begin{equation}
\frac{d^k G_{n}^{(\gamma)}}{dx^k} (x) = \frac{2^k \Gamma{(\gamma +k+1)}}
{\Gamma{(\gamma+1)}} G_{n-k}^{(\gamma +k)}(x),
\end{equation}
which coupled with (\ref{Gat1}) gives
\begin{equation}  \label{dkG1}
\frac{d^k G_{n}^{(\gamma)}}{dx^k} (1) = \frac{2^{k-1}
\Gamma{(\gamma +k)} \Gamma{(n+2\gamma
+k)}}{(n-k)!\Gamma{(\gamma+1)} \Gamma{(2\gamma +2k)}}, \qquad n
\ge 1.
\end{equation}
Gegenbauer polynomials correspond to Chebyshev polynomials of the 1st kind, 
$T_n(x)$, when $\gamma=0$, to Legendre $P_n(x)$ for $\gamma=1/2$
and to Chebyshev of the 2nd kind, $U_n(x)$, for $\gamma=1$. For
the non standard normalization,
\begin{equation} 
G_n^{(0)}(x) =\frac{T_n(x)}{n} , \quad G_n^{(1/2)}(x)=P_n(x),
\quad G_n^{(1)}(x)=\frac{U_n(x)}{2}.
\label{gegchebleg}
\end{equation}

\section{Integrals and Asymptotics}
\label{Asymptotics}
As shown in section \ref{Legendre}, the tau equations
(\ref{tauLeg}) provide a matrix eigenproblem of the form $\mu A a
= B a$. To reduce the coefficients $A(k,l)$ and $B(k,l)$ defined
in (\ref{Akl}) and (\ref{Bkl}) to the expressions (\ref{Akl2}) and
(\ref{Bkl2}),  use (\ref{gegsl}) and (\ref{gege2}) repeatedly
\begin{multline}\label{cl}
D^2 \left[ (1-x^2)^2 G_l^{(5/2)}(x)\right] = \frac{1}{5} D^2
\left[
(1-x^2)^2 D G_{l+1}^{(3/2)}(x)\right] \\
=-\frac{1}{5}(l+1)(l+4) D \left[ (1-x^2) G_{l+1}^{(3/2)}(x)\right] =
-\frac{1}{15}(l+1)(l+4) D \left[ (1-x^2) D G_{l+2}^{(1/2)}(x)\right] \\
= \frac{1}{15}(l+1)(l+2)(l+3) (l+4) \, G_{l+2}^{(1/2)}(x) \equiv \CC_l \, G_{l+2}^{(1/2)}(x).
\end{multline}
For the perturbation analysis described in section \ref{Legendre}
we need the first order corrections to $B(0,n-4)$ and to $B(1,n-5)$.
From equations (\ref{B1n-4}) and (\ref{B1n-5})
\begin{align}
B_1(0,n-4) = & \lim_{\epsilon \to 0}\; \frac{\CC_{n-4}}{\epsilon} \int_{-1}^1P_{n-2}(x)
\;(1-x^2)^{\epsilon} \;dx,  \label{B1}\\
B_1(1,n-5) = &\lim_{\epsilon \to 0} \; \frac{\CC_{n-5}}{\epsilon} \int_{-1}^1 x \; P_{n-3}(x)
\; (1-x^2)^{\epsilon} \;dx. \label{B2}
\end{align}
To evaluate $ \int_{-1}^1 P_n(x) (1-x^2)^{\epsilon}  dx$
 to $O(\epsilon)$ for $|\epsilon| \ll 1$ and $n$ even, use (\ref{gegsl}) for
$\gamma=1/2$ and integration by parts to derive
\begin{equation}
\begin{aligned}
\int_{-1}^1 P_n(x) (1-x^2)^{\epsilon}  dx = & \frac{-1}{n(n+1)}
\int_{-1}^1 D \left[ (1-x^2) D P_n\right]  (1-x^2)^{\epsilon} dx\\
=& \frac{\epsilon}{n(n+1)} \int_{-1}^1  DP_n \,
(1-x^2)^{\epsilon} (-2x) dx.
\end{aligned}
\end{equation}
This integral is 0 if $n$ is odd since $P_n(x)=(-1)^n P_n(-x)$.
Since we have an $\epsilon$ pre-factor we can now set $\epsilon=0$ in the integral, and
do the remaining integral by parts to obtain
\begin{equation}\label{Kn2}
\begin{aligned}
\int_{-1}^1 P_n(x) (1-x^2)^{\epsilon}  dx \sim& \frac{-2\epsilon}{n(n+1)} \int_{-1}^1  x DP_n dx  \\
&=  \frac{-2 \epsilon}{n(n+1)} \int_{-1}^1 \left( D(x P_n) - P_n \right)dx \\
&= \frac{-2 \epsilon}{n(n+1)} (P_n(1)+P_n(-1)) = \frac{-4
\epsilon}{n(n+1)}
\end{aligned}
\end{equation}
for $n$ even (0 for $n$ odd as should be).

For the integral in (\ref{B2}), use the recurrence (\ref{gegrec}) for $\gamma=1/2$ to write
$ (2n-5) x P_{n-3}(x) = (n-2) P_{n-2}(x) + (n-3) P_{n-4}(x)$ and evaluate the resulting 2 integrals from the formula (\ref{Kn2}).  Hence, for $n$ even,
\begin{align}
B_1(0,n-4) \sim & -\frac{4\CC_{n-4}}{(n-2)(n-1)},   \label{B1eval}\\
B_1(1,n-5) \sim & -\frac{4\CC_{n-5}}{(n-4)(n-1)}.  \label{B2eval}
\end{align}

Now for $A_0(0,n-4)$ and $A_0(1,n-5)$ and $n$ even, we have
\begin{align}
\label{An-4}
A_0(0,n-4)= & \CC_{n-4} \int_{-1}^1 D^2P_{n-2}(x) dx= (n-2)(n-1)\CC_{n-4}, \\
\label{An-5}
A_0(1,n-5)=& \CC_{n-5} \int_{-1}^1 xD^2P_{n-3}(x) dx=(n-4)(n-1) \CC_{n-5}.
\end{align}

\section{Proof of theorem \protect\ref{thm41}}

\label{ProofAppend}

Consider
\begin{equation}
f_n(x;z)=\sum_{k=0}^{\infty} z^k D^k G_n^{(\gamma)}(x) +
K\left(G_{n}^{(\gamma)}(x) - G_{n+2}^{(\gamma)}(x)\right)
\label{modeq50}
\end{equation}
where $z$ is a solution of $f_n(1;z)=0$ and $K$ is
\begin{equation}
K = \frac{(G_{n-1}^{(\gamma)}(1) -
G_{n+1}^{(\gamma)}(1))}{2(n+\gamma)(G_{n}^{(\gamma)}(1) -
G_{n+2}^{(\gamma)}(1))} = ... = \frac{n+2}{2(n+\gamma+1)(n+2\gamma-1)}
\end{equation}
where we have used (\ref{Gat1}). Note that
$f_n(1;z)=p^{(\gamma)}_n(z)$ defined in theorem \ref{thm41}.
Taking the $x$-derivative of $f_n(x;z)$ and using
(\ref{dgegrec}), we find
\begin{equation}
\frac{df_n(x;z)}{dx}=\sum_{k=0}^{\infty} z^k
D^{k+1}G_n^{(\gamma)}(x) - 2K(n+\gamma+1) G_{n+1}^{(\gamma)}(x).
\end{equation}
Thus $f_n(x;z)$ satisfies the following differential equation
\begin{equation}  \label{eqa10}
f_n(x;z)-(1+K)G_n^{(\gamma)}(x)+K G_{n+2}^{(\gamma)}(x) =
 z \frac{df_n(x;z)}{dx} +z 2K(n+\gamma+1) G_{n+1}^{(\gamma)}(x).
\end{equation}
Multiplying by $(1+x)\frac{d f^*_n(x;z)}{dx} $, integrating in
the Gegenbauer norm and adding the complex conjugate we get
\begin{multline}  \label{exp1}
\int_{-1}^1 \frac{d |f_n|^2}{dx} (1+x)w(x) dx
-(1+K)\left(\int_{-1}^1 (1+x)
\frac{df^*_n}{dx} G^{(\gamma)}_n(x) w(x) dx + C.C.\right) \\
+K\left(\int_{-1}^1 (1+x) \frac{df^*_n}{dx} G^{(\gamma)}_{n+2}(x)
w(x) dx + C.C.\right)= (z + z^*) \int_{-1}^1 \left|\frac{d
f_n}{dx}\right|^2
(1+x)w(x) dx \\
+\left( z 2K(n+\gamma+1)\int_{-1}^1 (1+x) \frac{d f^*_n}{dx}
G^{(\gamma)}_{n+1}(x)w(x) dx + C.C.\right)
\end{multline}
where $C.C.$ denotes the complex conjugate. To simplify
(\ref{exp1}) we need to compute four simple integrals
\begin{equation}
\begin{split}
I_1 = \int_{-1}^1 \frac{d |f_n|^2}{dx} (1+x)w(x) dx \\
J_0 = \int_{-1}^1 (1+x) \frac{df^*_n}{dx} G^{(\gamma)}_n(x) w(x) dx \\
J_1 = \int_{-1}^1 (1+x)\frac{d f^*_n}{dx} G^{(\gamma)}_{n+1}(x)w(x) dx \\
J_2 = \int_{-1}^1 (1+x) \frac{df^*_n}{dx} G^{(\gamma)}_{n+2}(x)
w(x) dx.
\end{split}
\end{equation}
Using integration by parts, the first integral becomes
\begin{equation}
I_1 = \int_{-1}^1 \frac{d |f_n|^2}{dx} (1+x)w(x) dx= -\int_{-1}^1
|f_n|^2(1-2\gamma x) \frac{w(x)}{1-x} dx.
\end{equation}
Therefore, the integral is negative for $-1/2 < \gamma \le 1/2$
since for this range of parameters $1-2\gamma x$ is positive.

For the calculation of the three other integrals we are going to
need the expression
\begin{align}
\frac{df_n}{dx}= & \sum_{k=0}^{\infty} z^k D^{k+1}
G^{(\gamma)}_{n}(x) -2K
(\gamma + n +1)G^{(\gamma)}_{n+1}(x) \\
=& \mathcal{P}_{n-2}(x;z) +2(\gamma+n-1) G^{(\gamma)}_{n-1}(x)
-2K (\gamma
+n+1) G^{(\gamma)}_{n+1}(x)  \label{fder2} \\
=& \mathcal{P}_{n-1}(x;z)-2K(\gamma+n+1) G^{(\gamma)}_{n+1}(x)
\label{fder}
\end{align}
where $\mathcal{P}_{n-2}(x;z)$ and $\mathcal{P}_{n-1}(x;z)$
are polynomials of degree $n-2$ and $n-1$, respectively.

With the use of (\ref{fder2}) and orthogonality of the Gegenbauer
polynomials, we find
\begin{multline}
J_0 = \int_{-1}^1 (1+x) \frac{df^*_n}{dx} G^{(\gamma)}_n(x) w(x) dx= \\
2(n-1+\gamma)\int_{-1}^1 x G^{(\gamma)}_{n-1}(x)G^{(\gamma)}_n(x)
w(x) dx -2K(n+1+\gamma)\int_{-1}^1
xG^{(\gamma)}_{n+1}(x)G^{(\gamma)}_{n}(x) w(x) dx
= \\
n \int_{-1}^1 (G^{(\gamma)}_n(x))^2 w(x) dx -
K(n+\gamma +1)\frac{n+1}
{n+\gamma}\int_{-1}^1 (G^{(\gamma)}_{n+1}(x))^2 w(x) dx= \\
n h_n^{(\gamma)}-K\frac{(n+\gamma +1)(n+1)}{n+\gamma} h_{n+1}^{(\gamma)} \\
\end{multline}
where we have also used (\ref{gege2}) and (\ref{gegrec}). In the
same way we compute
\begin{multline}
J_2 = \int_{-1}^1 (1+x) \frac{df^*_n}{dx} G^{(\gamma)}_{n+2}(x) w(x) dx= \\
-K \int_{-1}^1 2(n+\gamma +1) x G^{(\gamma)}_{n+1}(x)
G^{(\gamma)}_{n+2}(x) w(x) dx = -K(n+2)h_{n+2}^{(\gamma)}
\end{multline}
and
\begin{multline}
J_1 = \int_{-1}^1 (1+x)\frac{d f^*_n}{dx}
G^{(\gamma)}_{n+1}(x)w(x)
dx=-2K(n+\gamma +1)\int_{-1}^1 (1+x)(G^{(\gamma)}_{n+1}(x))^2 w(x) dx= \\
-2K(n+\gamma +1)\int_{-1}^1(G^{(\gamma)}_{n+1}(x))^2 w(x)
dx=-2K(n+\gamma +1)h_{n+1}^{(\gamma)}.
\end{multline}
Thus, (\ref{exp1}) transforms to
\begin{multline}  \label{exp2}
-\int_{-1}^1 |f_n|^2 \frac{(1-2\gamma x) w(x)}{(1-x)} dx
-2(1+K)\left(n h_{n}^{(\gamma)}-K\frac{(n+\gamma +1)(n+1)}{n+\gamma}
h_{n+1}^{(\gamma)} \right)\\ -2K^2(n+2)h_{n+2}^{(\gamma)} =(z+z^*)
\left( \int_{-1}^1 \left|\frac{d f_n}{dx}\right|^2 (1+x)w(x) dx -
4K^2 (n+1+\gamma)^2 h_{n+1}^{(\gamma)} \right).
\end{multline}
Our task is to show that the left hand side of the above
expression is negative whereas the coefficient of the term
$z+z^*$ on the right hand side is positive. A simple
calculation shows that
\begin{multline}
n h^{(\gamma)}_n-K \frac{(n+\gamma
+1)(n+1)}{(n+\gamma)}h^{(\gamma)}_{n+1}=
n h^{(\gamma)}_n- \frac{(n+1)(n+2)}{(n+\gamma)(n+2\gamma -1)}
h^{(\gamma)}_{n+1} \\
=\frac{\pi 2^{-1-2\gamma} \Gamma(n+2\gamma)}{\gamma^2
\Gamma^2(\gamma) n! (n+\gamma)}
\left(n-\frac{(n+2)(n+2\gamma)}{2(n-1+2\gamma)(n+1+\gamma)} \right)
\\ \ge \frac{\pi 2^{-1-2\gamma} \Gamma(n+2\gamma)}{\gamma^2
\Gamma^2(\gamma) n! (n+\gamma)}
\left(n-\frac{(n+2)(n+1)}{(n-2)(2n+1)} \right)\\
\end{multline}
The last parenthesis is positive for $n \ge 3$.

For the right hand side we use the notation in (\ref{fder}) to get
\begin{multline}
\int_{-1}^1 \left|\frac{d f_n}{dx}\right|^2 (1+x)w(x) dx - 4K^2
(n+1+\gamma)^2 h_{n+1}^{(\gamma)}= \\
\int_{-1}^1 (1+x)\left(\left|\mathcal{P}_{n-1}(x;z)\right|^2-
2K(n+1+\gamma)G^{(\gamma)}_{n+1}(x) \left(\mathcal{P}_{n-1}(x;z)
+\mathcal{P}_{n-1}^*(x;z) \right)\right. \\
\left. +4K^2(n+1+\gamma)^2 (G^{(\gamma)}_{n+1}(x))^2 \right) w(x)
dx- 4K^2
(n+1+\gamma)^2 h_{n+1}^{(\gamma)}= \\
\int_{-1}^1(1+x)|\mathcal{P}_{n-1}(x;z)|^2 w(x) dx +0+0+4K^2
(n+1+\gamma)^2 h_{n+1}^{(\gamma)}- 4K^2 (n+1+\gamma)^2 h_{n+1}^{(\gamma)} \\
=\int_{-1}^1(1+x)\left|\mathcal{P}_{n-1}(x;z)\right|^2 w(x) dx >
0.
\end{multline}
Thus (\ref{exp2}) becomes
\begin{multline}  \label{exp3}
-\int_{-1}^1 |f_n|^2 \frac{(1-2\gamma x) w(x)}{(1-x)} dx
-2(1+K)\left(n h_{n}^{(\gamma)}- K\frac{(n+\gamma
+1)(n+1)}{n+\gamma} h_{n+1}^{(\gamma)}\right)\\
-2K^2(n+2)h_{n+2}^{(\gamma)}
=(z+z^*)\int_{-1}^1(1+x)\left|\mathcal{P}_{n-1}(x;z)\right|^2 w(x)
dx
\end{multline}
and $\Re(z)<0$.

For $n=2$ we get the following characteristic polynomial
\begin{equation}
f^{(\gamma)}_2(1;z) = \frac{2}{3}(\gamma+1)(3z^2 + 3z + 1)
\end{equation}
whose zeros
\begin{equation}
z_1 = -\frac{1}{2} - \frac{\sqrt{3}i}{6}, \quad \quad z_2
=-\frac{1}{2} + \frac{\sqrt{3}i}{6}
\end{equation}
have negative real parts for any $\gamma$.

\bibliographystyle{siam}
\bibliography{../MariosThesis}

\end{document}